\documentclass[12pt]{article}
\usepackage{amssymb}
\usepackage{amsfonts}
\usepackage{latexsym}
\usepackage[dvips]{graphics}

\newtheorem{thm}{Theorem}[section]
\newtheorem{lem}[thm]{Lemma}
\newtheorem{alg}[thm]{Algorithm}
\newtheorem{cor}[thm]{Corollary}

\newtheorem{rem}[thm]{Remark}

\newtheorem{ex}[thm]{Example}

\large\normalsize

\begin{document}

\begin{center}
{\Large\bf Enumeration of subtrees of trees}
\\[30pt]
{Weigen\ Yan$^{\rm a,b}$ \footnote{This work is supported by
FMSTF(2004J024) and NSFF(E0540007).} \quad and \quad Yeong-Nan\
Yeh$^{\rm b}$ \footnote{Partially supported by NSC
95-2115-M-001-009.\\
\hspace*{5mm} {\it Email address:} weigenyan@263.net (W. G. Yan),
mayeh@math.sinica.edu.tw (Y. N. Yeh).}}
\\[10pt]
\footnotesize { $^{\rm a}$School of Sciences, Jimei University,
Xiamen 361021, China
\\[7pt]
$^{\rm b}$Institute of Mathematics, Academia Sinica, Taipei 11529.
Taiwan.}
\end{center}
\begin{abstract}
Let $T$ be a weighted tree. The weight of a subtree $T_1$ of $T$
is defined as the product of weights of vertices and edges of
$T_1$. We obtain a linear-time algorithm to count the sum of
weights of subtrees of $T$. As applications, we characterize the
tree with the diameter at least $d$, which has the maximum number
of subtrees, and we characterize the tree with the maximum degree
at least $\Delta$, which has the minimum number of subtrees.
\\
\\
{\sl Keywords:}\quad subtree, extremal tree, tree transformation,
diameter, connected subgraph
\end{abstract}
\section{Introduction}
\hspace*{\parindent}
Throughout this paper, we suppose that $T=(V(T),E(T);f,g)$ is a
weighted tree with the vertex set $V(T)=\{v_1,v_2,\cdots,v_n\}$,
the edge set $E(T)=\{e_1,e_2,\cdots, e_{n-1}\}$, vertex-weight
function $f: V(T)\rightarrow \mathcal R$ and edge-weight function
$g: E(T)\rightarrow \mathcal R$ (where $\mathcal R$ is a
commutative ring with a unit element $1$), if not otherwise
specified. If a weighted tree $T=(V(T),E(T);f,g)$ satisfies
$f=g=1$, we call $T$ a simple tree and denote it by
$T=(V(T),E(T))$. Let $\mathcal T(T)$ denote the set of subtrees of
a tree $T$. For arbitrary two fixed vertices $v_i$ and $v_j$,
denote by $\mathcal T(T;v_i)$ (resp. $\mathcal T(T;v_i,v_j)$) the
set of subtrees of $T$, each of which contains vertex $v_i$ (resp.
vertices $v_i$ and $v_j$), denote by $a(T;k)$ the number of
subtrees of $T$ with $k$ edges, denote by $a(T;v_i;k)$ (resp.
$a(T;v_i,v_j;k)$) the number of subtrees of $T$, each of which
contains vertex $v_i$ (resp. vertices $v_i$ and $v_j$) and $k$
edges, denote by $b(T;k)$ the number of subtrees of $T$ with $k$
vertices, and denote by $b(T;v_i;k)$ (resp. $b(T;v_i,v_j;k)$) the
number of subtrees of $T$ with $k$ vertices, each of which
contains vertex $v_i$ (resp. vertices $v_i$ and $v_j$).
Obviously, for any $k=0,1,\ldots, n-1$, we have the following:
$$a(T;k)=b(T;k+1), a(T;v_i;k)=b(T;v_i;k+1), a(T;v_i,v_j;k)=b(T;v_i,v_j;k+1). $$
For a given subtree $T_1$ of a weighted $T$, we define the weight
of $T_1$, denoted by $\omega(T_1)$, as the product of the weights
of the vertices and edges in $T_1$. The generating function of
subtrees of a weighted tree $T=(V(T),E(T);f,g)$, denoted by
$F(T;f,g)$, is the sum of weights of subtrees of $T$. That is,
$F(T;f,g)=\sum\limits_{T_1\in \mathcal T(T)}\omega(T_1).$
Similarly, we can define the generating function of subtrees of a
weighted tree $T=(V(T),E(T);f,g)$ containing a fixed vertex $v_i$
(resp. two fixed vertices $v_i$ and $v_j$), as the sum of weights
of subtrees of $T$ containing vertex $v_i$ (resp. vertices $v_i$
and $v_j$), denoted by $F(T;f,g;v_i)$ (resp. $F(T;f,g;v_i,v_j)$).
Hence we have
$$F(T;f,g;v_i)=\sum_{T_1\in \mathcal T(T;v_i)}\omega(T_1),\ \
\ \ F(T;f,g;v_i,v_j)=\sum_{T_1\in \mathcal
T(T;v_i,v_j)}\omega(T_1).
$$

By the definitions of $F(T;f,g), F(T;f,g;v_i)$ and
$F(T;f,g;v_i,v_j)$, if we weight each edge by $x$ and each vertex
by $y$, then
$$
\begin{array}{ll}
&F(T;y,x)=\sum\limits_{k=0}^{n-1}a(T;k)x^ky^{k+1}=\sum\limits_{k=1}^{n}b(T;k)x^{k-1}y^{k};\\
&F(T;y,x;v_i)=\sum\limits_{k=0}^{n-1}a(T;v_i;k)x^ky^{k+1}=\sum\limits_{k=1}^{n}b(T;v_i;k)x^{k-1}y^{k};\\
&F(T;y,x;v_i,v_j)=\sum\limits_{k=0}^{n-1}a(T;v_i,v_j;k)x^ky^{k+1}=\sum\limits_{k=1}^{n}b(T;v_i,v_j;k)x^{k-1}y^{k}.
\end{array}
$$

Let $T$ be a simple tree of order $n$, and let $v_i$ and $v_j$ be
arbitrary two distinct vertices of $T$. For the sake of
convenience, we denote by $\chi(T)=F(T;1,1)$ the number of
subtrees of $T$, by $\chi(T;v_i)=F(T;1,1;v_i)$ the number of
subtrees of $T$, each of which contains vertex $v_i$, and by
$\chi(T;v_i,v_j)=F(T;1,1;v_i,v_j)$ the number of subtrees of $T$,
each of which contains vertices $v_i$ and $v_j$.

Sz\'ekely and Wang \cite{SW05} studied the problem enumerating
subtrees of a tree. They proved the following:
\begin{thm}[Sz\'ekely and Wang \cite{SW05}]
The path $P_n$ has ${n+1\choose 2}$ subtrees, fewer than any other
trees of $n$ vertices. The star $K_{1,n-1}$ has $2^{n-1}+n-1$
subtrees, more than any other trees of $n$ vertices.
\end{thm}

Sz\'ekely and Wang \cite{SW05} said that it was not difficult to
design a recursive algorithm that would compute the number of
subtrees of a tree in a time bounded by a polynomial of $n$, the
number of vertices (but we have not found such an algorithm).
These may be the first results on enumeration of subtrees of a
simple tree. For some related results see also Sz\'ekely and Wang
\cite{SW052,SW053} and Wang \cite{W05}.

In the next section, we give a linear-time algorithm to count the
generating functions $F(T;f,g),F(T;f,g;v_i)$, and
$F(T;f,g;v_i,v_j)$ of subtrees of a weighted tree
$T=(V(T),E(T);f,g)$ for any two vertices $v_i$ and $v_j$. As an
applications, in Section 3 we characterize the tree with the
diameter at least $d$, which has the maximum number of subtrees,
and we characterize the tree with the maximum degree at least
$\Delta$, which has the minimum number of subtrees. Finally,
Section 4 presents our conclusions.
\section{Algorithms}
\hspace*{\parindent}
Let $T=(V(T),E(T);f,g)$ be a weighted tree of order $n>1$ and $u$
a pendant vertex of $T$. Suppose $e=(u,v)$ is the pendant edge of
$T$. We define a weighted tree $T'=(V(T'),E(T');f',g')$ of order
$n-1$ from $T$ as follows: $V(T')=V(T)\backslash \{u\},\
E(T')=E(T)\backslash \{e\}$, and
$$f'(v_s)=\left \{
\begin{array}{ll}
f(v)(f(u)g(e)+1) & \ \mbox{if}\ v_s=v,\\
f(v_s) & \ \mbox{otherwise}.
\end{array}
\right.,$$ for any $v_s\in V(T')$, and $g'(e)=g(e)$ for any $ e\in
E(T')$. Figure 1 illustrates the procedure constructing $T'$ from
$T$.
\begin{figure}[htbp]
  \centering
 \includegraphics{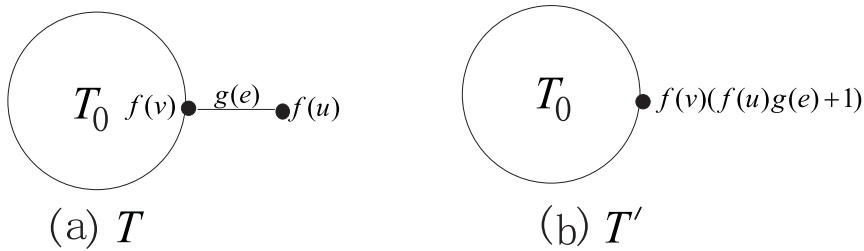}
  \caption{\ (a)\ A weighted tree $T=(V(T),E(T);f,g)$
with a pendent edge $e=(u,v)$.\ (b)\ The corresponding weighted
tree $T'=(V(T'),E(T');f',g')$.}
\end{figure}
\begin{thm}
Keeping the above notation, we have
$$F(T;f,g)=F(T';f',g')+f(u). \eqno{(1)}$$
\end{thm}

{\bf Proof}\ \ We partition the sets $\mathcal T(T)$ and $\mathcal
T(T')$ of subtrees of $T$ and $T'$ as follows:
$$\mathcal T(T)=\mathcal T_1\cup \mathcal T_{1'}\cup \mathcal T_2\cup
\mathcal T_3,\ \ \mathcal T(T')=\mathcal T_1'\cup \mathcal T_2',\
\mbox{where}$$

$\mathcal T_1$ is the set of subtrees of $T$, each of which
contains vertex $v$ but not vertex $u$;

$\mathcal T_{1'}$ is the set of subtrees of $T$, each of which
contains edges $e=(u,v)$;

$\mathcal T_2$ is the set of subtrees of $T$, each of which
contains neither $u$ nor $v$;

$\mathcal T_3$ is the set of subtrees of $T$, each of which
contains $u$ but not $v$;

$\mathcal T_1'$ is the set of subtrees of $T'$, each of which
contains vertex $v$;

$\mathcal T_2'$ is the set of subtrees of $T'$, each of which
contains no vertex $v$.

By the definitions above, we have

$(i)$\ \ there exist two natural bijections (ignore weights)
$\theta_1: T_1\longmapsto T_1'$ between $\mathcal T_1$ and
$\mathcal T_1'$, and $\theta_2: T_2\longmapsto T_2'$ between
$\mathcal T_2$ and $\mathcal T_2'$;

$(ii)$\ \ $\mathcal T_{1'}=\{T_1+u|T_1\in \mathcal T_1\}$, where
$T_1+u$ is the tree obtained from $T_1$ by attaching a pendant
edge $(v,u)$ at vertex $v$ of $T_1$;

$(iii)$\ \ $\mathcal T_3=\{u\}$.

Note that we have
$$\sum_{T_1'\in \mathcal T_1'}\omega(T_1')=\sum_{T_1'\in
\mathcal T_1'}f'(v)\frac{\omega(T_1')}{f'(v)}=\sum_{T_1'\in
\mathcal T_1'}f(v)[f(u)g(e)+1]\frac{\omega(T_1')}{f'(v)}.
\eqno{(2)}$$

By $(i),(ii)$ and $(iii)$, we have
$$
\sum_{T_{1'}\in \mathcal T_{1'}} \omega(T_{1'})=\sum_{T_1\in
\mathcal T_1}f(u)g(e)\omega(T_1), \eqno{(3)}$$
$$\sum_{T_2'\in \mathcal T_2'}\omega(T_2')=\sum_{T_2\in
\mathcal T_2}\omega(T_2),\eqno{(4)}$$
$$\sum_{T_3\in \mathcal T_3}\omega(T_3)=f(u).\eqno{(5)}$$
By $(3)$, we have
$$\sum_{T_1\in \mathcal T_1}\omega(T_1)+\sum_{T_{1'}\in
\mathcal T_{1'}}\omega(T_{1'})=\sum_{T_1\in \mathcal
T_1}[f(u)g(e)+1] \omega(T_1)=\sum_{T_1\in \mathcal
T_1}f(v)[f(u)g(e)+1] \frac{\omega(T_1)}{f(v)}. \eqno{(6)}$$ By
$(i)$, $\theta_1: T_1\longmapsto T_1'$ is a natural bijection
between $\mathcal T_1$ and $\mathcal T_1'$, then
$\frac{\omega(T_1')}{f'(v)}=\frac{\omega(T_1)}{f(v)}$ since $T_1$
and $T_1'$ have ``almost all" the same weights of vertices and
edges except the weights of $v$ in $T_1$ and $T_1'$ (one is $f(v)$
and another is $f(v)(f(u)g(e)+1)$). So by $(2)$ and $(6)$ we have
$$\sum_{T_1\in \mathcal T_1}\omega(T_1)+\sum_{T_{1'}\in
\mathcal T_{1'}}\omega(T_{1'})=\sum_{T_1'\in \mathcal
T_1'}\omega(T_1').\eqno{(7)}$$ Hence by $(4),(5),(7)$, and the
definitions of $F(T;f,g)$ and $F(T';f',g')$ we have
$$
F(T;f,g)=\sum\limits_{T_1\in \mathcal
T_1}\omega(T_1)+\sum\limits_{T_{1'}\in \mathcal
T_{1'}}\omega(T_{1'})+\sum\limits_{T_2\in \mathcal
T_2}\omega(T_2)+\sum\limits_{T_3\in \mathcal T_3}\omega(T_3)$$
$$
=\sum\limits_{T_1'\in \mathcal
T_1'}\omega(T_1')+\sum\limits_{T_2'\in \mathcal
T_2'}\omega(T_2')+f(u)=F(T',f',g')+f(u),
$$
and the theorem thus follows.\ \ \ \ $\blacksquare$

By a similar argument we have the following:
\begin{thm}
Let $T=(V(T),E(T);f,g)$ be a weighted tree of order $n>1$ and $u$
a pendant vertex of $T$. Suppose $e=(u,v)$ is the pendant edge of
$T$. Let $T'$ be the weighted tree defined as above. Then, for
arbitrary vertex $v_i\neq u$,  the generating functions
$F(T;f,g;v_i)$  and $F(T';f',g';v_i)$ of subtrees of $T$ and $T'$
satisfy the following:
$$F(T;f,g;v_i)=F(T';f',g';v_i). \eqno{(8)}$$
\end{thm}
\begin{thm}
Let $T=(V(T),E(T);f,g)$ be a weighted tree of order $n>1$ and $u$
a pendant vertex of $T$. Suppose $e=(u,v)$ is the pendant edge of
$T$. Let $T'$ be the weighted tree defined as above. Then, for
arbitrary two distinct vertices $v_i$ and $v_j$ such that $v_i\neq
u,v_j\neq u$,  the generating functions $F(T;f,g;v_i,v_j)$ and
$F(T';f',g';v_i,v_j)$ of subtrees of $T$ and $T'$ satisfy the
following:
$$F(T;f,g;v_i,v_j)=F(T';f',g';v_i,v_j). \eqno{(9)}$$
\end{thm}

For the sake of convenience, if $\{a_n\}_{\geq 0}$ is a sequence,
we define:  $\prod\limits_{t=i}^j a_t=1$ if $j<i$.
\begin{cor}
Let $P_n=(V(P_n),E(P_n);f,g)$ be a weighted path of order $n$,
where $V(P_n)=\{v_i|i=1,2,\ldots,n\},
E(P_n)=\{e_i=(v_i,v_i+1)|i=1,2,\ldots,n-1\}$, $f(v_i)=y_i$ for
$i=1,2,\ldots,n$, and $g(e_i)=x_i$ for $i=1,2,\cdots,n-1$. Then
$$F(P_n;f,g)=\sum_{j=0}^{n-1}\sum_{i=1}^{n-j}\left(\prod_{s=i}
^{i+j-1}x_sy_s\right)y_{i+j}, \ \eqno{(10)}$$
$$F(P_n;f,g;v_1)=y_1[1+\sum_{j=1}^{n-1}
\prod_{i=1}^j(x_iy_{i+1})]. \eqno{(11)}$$
\end{cor}

{\bf Proof}\ \ We prove the corollary by induction on $n$. It is
easy to prove that if $n=2$ or $3$ the corollary holds. Now we
suppose $n>3$ and proceed by induction. Let
$P_{n-1}'=(V(P_{n-1}'),E(P_{n-1}');f',g')$, where
$V(P_{n-1}')=\{v_i|i=1,2,\ldots,n-1\},
E(P_{n-1}')=\{e_i=(v_i,v_{i+1})|i=1,2,\ldots,n-2\}$, $f'(v_i)=y_i$
for $i=1,2,\ldots,n-2$ and $f'(v_{n-1})=y_{n-1}(y_{n}x_{n-1}+1)$,
and $g'(e_i)=x_i$ for $i=1,2,\cdots,n-2$. Then, by Theorem 2.1, we
have
$$F(P_n;f,g)=F(P_{n-1}';f',g')+y_n. $$
By induction, we have
$$F(P_{n-1}';f',g')=\sum_{j=0}^{n-2}\sum_{i=1}^{n-1-j}
\left(\prod_{s=i}^{i+j-1}x_sy_s'\right)y_{i+j}',$$ where
$y_s'=y_s$ for $s=1,2,\ldots,n-2$, and
$y_{n-1}'=y_{n-1}(y_{n}x_{n-1}+1)$. Hence  we have
$$
F(P_n;f,g)=F(P_{n-1}';f',g')+y_n=\sum_{j=0}^{n-2}
\sum_{i=1}^{n-1-j}\left(\prod_{s=i}^{i+j-1}x_sy_s'\right)y_{i+j}'+y_n
$$
$$
=\sum\limits_{i=1}^{n-1}y_i'+\sum\limits_{i=1}^{n-2}y_i'y_{i+1}'x_i
+\ldots+\sum\limits_{i=1}^{n-1-k}y_i'y_{i+1}'\ldots
y_{i+k}'x_ix_{i+1}\ldots x_{i+k-1}$$$$+\ldots+y_1'y_2'\ldots
y_{n-1}'x_1x_2\ldots x_{n-2}+y_n. \eqno{(12)}$$ Note that
$y_i'=y_i$ for $i=1,2,\ldots,n-2$, and
$y_{n-1}'=y_{n-1}(y_{n}x_{n-1}+1)$. By $(12)$, it is easy to show
that $(10)$ holds. Similarly, we can show that $(11)$ holds and
hence the corollary has been proved.\ \ \ \ $\blacksquare$

A direct result of Corollary 2.4 is the following:
\begin{cor}
$$F(P_n;y,x)=\sum_{j=0}^{n-1}(n-j)x^{j}y^{j+1},\ \
F(P_n;y,x;v_1)=\sum_{j=0}^{n-1}x^jy^{j+1},$$
$$F(P_n;y,1)=\sum_{j=1}^{n}(n-j+1)y^j,\ F(P_n;1,x)
=\sum_{j=0}^n(n-j)x^j.$$
\end{cor}

Similarly, we can prove the following:
\begin{cor}
Let $K_{1,n-1}=(V(K_{1,n-1}),E(K_{1,n-1});f,g)$ be a weighted star
of order $n$, where $V(K_{1,n-1})=\{v_i|i=1,2,\ldots,n\},
E(K_{1,n-1})=\{e_i=(v_n,v_{i})|i=1,2,\ldots,n-1\}$, $f(v_i)=y_i$
for $i=1,2,\ldots,n$, and $g(e_i)=x_i$ for $i=1,2,\cdots,n-1$.
Then
$$F(K_{1,n-1};f,g)=\sum_{i=1}^n y_i+\sum_{i=1}^{n-1}\left[
\sum_{1\leq j_1<j_2<\ldots<j_i\leq n-1}\left(\prod_{k=1}^{i}
x_{j_k}y_{j_k}\right)\right ]y_n.$$
\end{cor}
\begin{cor}
$$F(K_{1,n-1};y,x)=ny+\sum_{i=1}^{n-1}{n-1\choose i}x^{i}y^{i+1}. $$
\end{cor}

By Corollaries 2.5 and 2.7, we have the following:
\begin{cor} [Sz\'ekely and Wang \cite{SW05}]
$$\chi(P_n)=F(P_n;1,1)={n+1\choose 2},\ \ \chi(K_{1,n-1})=F(K_{1,n-1};1,1)=2^{n-1}+n-1.$$
\end{cor}

By Theorems 2.1, 2.2, and 2.3, we can produce three
graph-theoretical algorithms for computing the generating
functions $F(T;f,g), F(T;f,g;v_i)$, and $F(T;f,g;v_i,v_j)$ of
subtrees of a weighted tree $T=(V(T),E(T);f,g)$ directly from $T$
for arbitrary two different vertices $v_i$ and $v_j$,
respectively, as follows:
\begin{alg}
Let $T=(V(T),E(T);f,g)$ be a weighted tree with two or more
vertices.

{\bf Step 1}\ \ Initialize.

Define: $p(v_s)=f(v_s)$, for all $v_s\in V(T)$; and $N=0$.

{\bf Step 2}\ \ Contract.

(a)\ \ Choose a pendant vertex $u$ and let $e=(u,v)$ denote the
pendant edge.

(b)\ \ Replace $p(v)$ with $p(v)(p(u)g(e)+1)$.

(c)\ \ Replace $N$ with $N+p(u)$.

(d)\ \ Eliminate vertex $u$ and edge $e$.

{\bf Step 3}\ \ If $v$ is the only remaining vertex, go to Step 4.
Otherwise, go to Step 2.

{\bf Step 4}\ \ Answer: $F(T;f,g)=p(v)+N$.
\end{alg}
\begin{alg}
Let $T=(V(T),E(T);f,g)$ be a weighted tree with two or more
vertices and $v_i$ a fixed vertex of $T$.

{\bf Step 1}\ \ Initialize.

Define: $p(v_s)=f(v_s)$, for all $v_s\in V(T)$.

{\bf Step 2}\ \ Contract.

(a)\ \ Choose a pendant vertex $u\neq v_i$ and let $e=(u,v)$
denote the pendant edge.

(b)\ \ Replace $p(v)$ with $p(v)(p(u)g(e)+1)$.

(c)\ \ Eliminate vertex $u$ and edge $e$.

{\bf Step 3}\ \ If $v$ is the only remaining vertex $v_i$, go to
Step 4. Otherwise, go to Step 2.

{\bf Step 4}\ \ Answer: $F(T;f,g;v_i)=p(v)$.
\end{alg}
\begin{alg}
Let $T=(V(T),E(T);f,g)$ be a weighted tree with two or more
vertices, and $v_i$ and $v_j$ two distinct vertices of $T$.

{\bf Step 1}.\ \ Initialize.

\ \ \ \ Define: $p(v_s)=f(v_s),$ for all $v_s\in V(T)$.

{\bf Step 2}\ \ If $T$ is a path, and $v_i$ and $v_j$ are two
pendant vertices, go to Step 5. Otherwise, go to Step 3.

{\bf Step 3}\ \ Contract.

$(a)$\ \ Choose a pendant vertex $u$, which is different from
$v_i$ and $v_j$, and let $e=(u,v)$ denote the pendant edge.

$(b)$\ \ Replace $p(v)$ with $p(v)(p(u)g(e)+1)$.

$(c)$\ \ Eliminate vertex $u$ and edge $e$.

{\bf Step 4}\ \ If there exists no vertex $u$ satisfying the
condition $(a)$ in Step 3, go to Step 5. Otherwise, go to Step 3.

{\bf Step 5}\ \ Answer: $F(T;f,g;v_i,v_j)=\prod\limits_{v\in
V(P_{v_iv_j})}p(v)\prod\limits_{e\in E(P_{v_iv_j})}g(e)$, where
$P_{v_iv_j}$ denotes the unique path of $T$ from vertex $v_i$ to
$v_j$.
\end{alg}
\begin{rem}
It is not difficult to see that Algorithms 2.9, 2.10, and 2.11 are
linear on the number of vertices of the tree $T$. Let $T$ be a
simple tree of order $n$ and $v_i$ and $v_j$ two distinct vertices
of $T$. By Algorithms 2.9, 2.10, and 2.11, we can compute easily
the numbers $\chi(T), \chi(T;v_i), \chi(T;v_i,v_j), a(T;k),
b(T;k), a(T;v_i;k), a(T;v_i,v_j;k), b(T;v_i;k)$ and
$b(T;v_i,v_j;k)$, respectively. The following examples show these
procedures of computations.
\end{rem}
\begin{ex}
We compute the numbers $\chi(T), \chi(T;B), \chi(T;A,B)$ of a
simple tree $T$, which appears in the upper left corner in Figure
2. We weight each vertex and edge of $T$ by one. From the
illustration in Figure 2, we know that $\chi(T)=62,
\chi(T;B)=24(1\times 1+1)=48, \chi(T;A)=25, \chi(T;A,B)=1\times
1\times 24=24$.
\end{ex}
\begin{figure}[htbp]
  \centering
 \scalebox{0.9}{\includegraphics{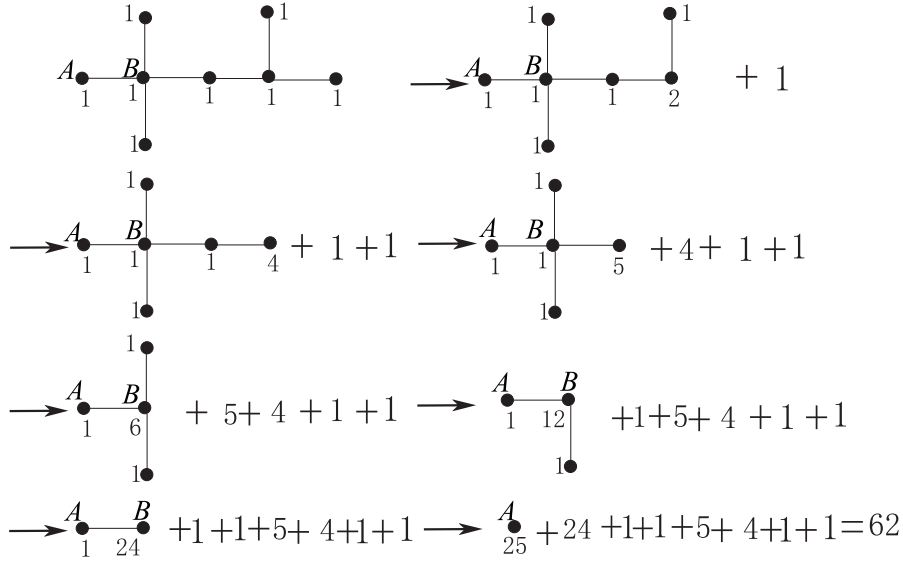}}
  \caption{\ An illustration of the procedures for
computing the numbers $\chi(T), \chi(T;B), \chi(T;A,B)$ of a
simple tree by Algorithms 2.9, 2.10, and 2.11.}
\end{figure}
\begin{ex}
We compute the edge generating functions $F(T;1,x), F(T;1,x;A)$
and $F(T;1,x;B,C)$ of a simple tree $T$, which appears in Figure
3. We can weight each vertex by one and each edge by $x$ (or
weight each vertex by $y$ and each edge by one, see Example 2.15).
From the illustration in Figure 3, we know that
$F(T;1,x)=x(x^2+2x+1)^2+2(x^2+2x+1)+4=x^5+4x^4+6x^3+6x^2+5x+6,
F(T;1,x;A)=x(x^2+2x+1)^2+(x^2+2x+1)=x^5+4x^4+6x^3+5x^2+3x+1,
F(T;1,x;B,C)=x(x+1)x(x^2+2x+1)=x^5+3x^4+3x^3+x^2$. Hence
$a(T;0)=6,a(T;1)=5,a(T;2)=6,a(T;3)=6,a(T;4)=4,a(T;5)=1;
a(T;A;0)=1,a(T;A;1)=3,a(T;A;2)=5,a(T;A;3)=6,a(T;A;4)=4,a(T;A;5)=1;
a(T;B,C;0)=0,
a(T;B,C;1)=0,a(T;B,C;2)=1,a(T;B,C;3)=3,a(T;B,C;4)=3,a(T;B,C;5)=1$.
\end{ex}
\begin{figure}[htbp]
  \centering
 \scalebox{1.1}{\includegraphics{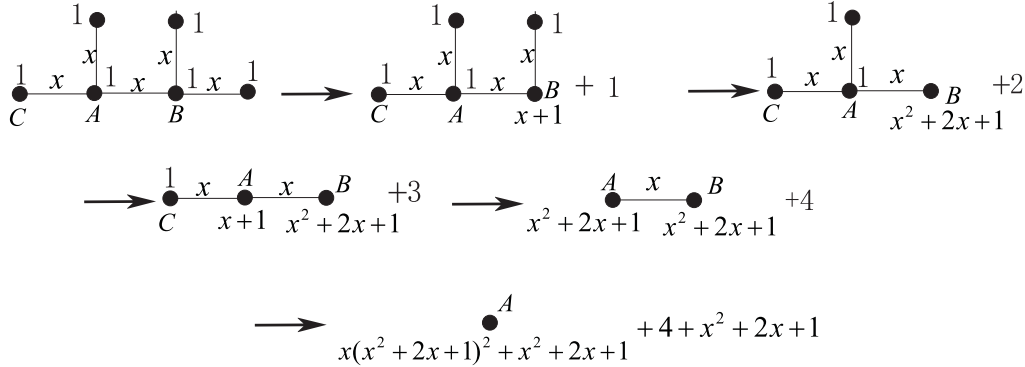}}
  \caption{\ An illustration of the procedures for
computing $F(T;1,x), F(T;1,x;A)$ and $F(T;1,x;B,C)$ of a simple
tree $T$ by Algorithms 2.9, 2.10, and 2.11.}
\end{figure}
\begin{ex}
We compute the vertex generating functions $F(T;y,1), F(T;y,1;A)$
and $F(T;y,1;B,C)$ of a simple tree $T$, which appears in Figure 4
or Figure 3. We weight each vertex by $y$ and each edge by $1$.
From the illustration in Figure 4, we know that
$F(T;y,1)=(y^3+2y^2+y)^2+2(y^3+2y^2+y)+4y=y^6+4y^5+6y^4+6y^3+5y^2+6y,
F(T;y,1;A)=(y^3+2y^2+y)^2+y^3+2y^2+y=y^6+4y^5+6y^4+5y^3+3y^2+y,
F(T;y,1;B,C)=y(y^2+y)(y^3+2y^2+y)=y^6+3y^5+3y^4+y^3$. Hence
$b(T;1)=6,b(T;2)=5,b(T;3)=6,b(T;4)=6,b(T;5)=4,b(T;6)=1;
b(T;A;1)=1,b(T;A;2)=3,b(T;A;3)=5,b(T;A;4)=6,b(T;A;5)=4,b(T;A;6)=1;
b(T;B,C;1)=0,b(T;B,C;2)=0,b(T;B,C;3)=1,b(T;B,C;4)=3,b(T;B,C;5)=3,b(T;B,C;6)=1$.
\end{ex}
\begin{figure}[htbp]
  \centering
 \scalebox{1.1}{\includegraphics{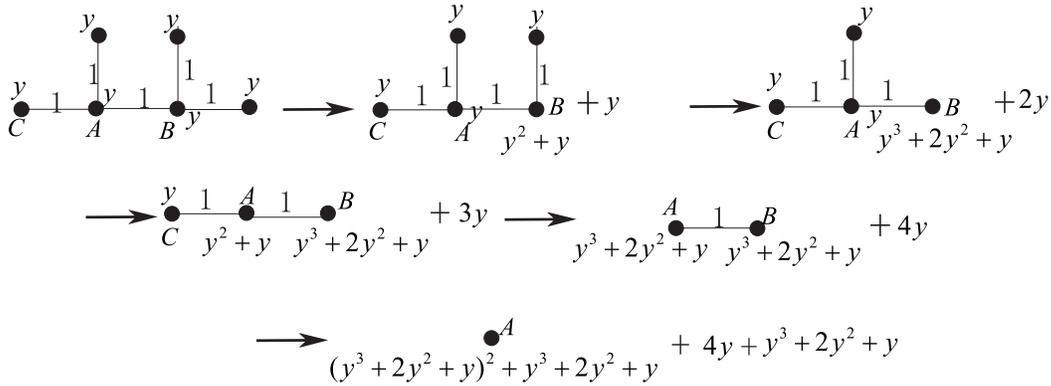}}
  \caption{\ An illustration of the procedures for
computing $F(T;y,1), F(T;y,1,;A)$ and $F(T;y,1;B,C)$ of a simple
tree $T$ by Algorithms 2.9, 2.10, and 2.11.}
\end{figure}

From Example 2.14, for the tree $T$ shown in Figure 3, we have
$\chi(T)=\sum\limits_{k=0}^5a(T;k)=28,
\chi(T;A)=\sum\limits_{k=0}^5a(T;A;k)=20,\chi(T;B,C)=\sum\limits_{k=0}^5a(T;B,C;k)=8$.
\section{Trees with extremal number of subtrees}
\hspace*{\parindent}
We suppose that the tree $T$ considered in this section is simple,
if not specified. In Section 3.1, we introduce four
transformations of trees, each of which gives us a way of
comparing numbers of subtrees of a pair of trees. In Section 3.2,
by the four transformations of trees we characterize the tree with
the diameter at least $d$, which has the maximum number of
subtrees, and we also characterize the tree with the maximum
degree at least $\Delta$, which has the minimum number of
subtrees. As corollaries, we obtain the trees with the second,
third, fourth, and fifth largest numbers of subtrees and the tree
with the second minimum number of subtrees.
\subsection{Four transformations of trees} \hspace*{\parindent}
Denote the degree of a vertex $v$ of tree $T$ by $d_T(v)$. Let
$T_1'$ and $T_2'$ be two trees, and let $u$ (resp. $v$) be a
vertex of $T_1'$ (resp. $T_2'$), where $|V(T_2')|=r+1\geq 2$. Let
$T_1$ be a tree obtained from $T_1'$ and $T_2'$ by identifying
vertices $u$ and $v$ (see the illustration in Figure 5(a)).
Construct a tree $T_2$ from $T_1'$ by attaching $r$ pendant edges
to vertex $u$ of $T_1'$ (see Figure 5(b)). We call the procedure
constructing $T_2$ from $T_1$ the first transformation of tree
$T_1$, denoted by $\phi_{1}(T_1)=T_2$.
\begin{figure}[htbp]
  \centering
 \scalebox{0.8}{\includegraphics{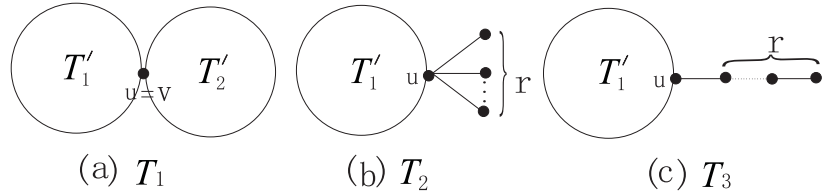}}
  \caption{\ \ (a)\ The tree $T_1$. \ (b)\ The tree
$T_2$.\ (c)\ The tree $T_3$.}
\end{figure}
\begin{lem}
Let $T_1$ and $T_2$ be the trees defined as above, where $r\geq 1$
and $|V(T_1')|\geq 2$ . Then
$$\chi(T_1)=F(T_1;1,1)\leq \chi(T_2)=F(T_2;1,1)$$
with equality holds if and only if $T_2'=K_{1,r}$ and
$d_{T_2'}(v)=r$.
\end{lem}

{\bf Proof}\ \ Let $f_i: V(T_1')\longrightarrow \mathcal R$
($i=1,2$) be two functions defined as follows:
$$
f_1(v')=\left\{
\begin{array}{ll}
F(T_2';1,1;v) & \ \mbox{if}\ \ v'=u,\\
1 & \ \mbox{otherwise}.
\end{array}
\right. , f_2(v')=\left\{
\begin{array}{ll}
2^r & \ \mbox{if}\ \ v'=u,\\
1 & \ \mbox{otherwise}.
\end{array}
\right.,
$$
where $F(T;1,1;v)$ is the number of subtrees of $T$, each of which
contains vertex $v$. Suppose that $\Phi_u(T_1')$ is the set of
subtrees of $T_1'$ with as least two vertices, each of which
contains vertex $u$. By Algorithms 2.9 and 2.10, we have
$$
\begin{array}{lll}
F(T_1;1,1)&=&F(T_2';1,1)-F(T_2';1,1;v)+F(T_1';f_1,1)\\
&=&F(T_2';1,1)-F(T_2';1,1;v)+F(T_1'-u;1,1)+F(T_2';1,1;v)[1+|\Phi_u(T_1')|]\\
&=&F(T_2';1,1)+F(T_1'-u;1,1)+F(T_2';1,1;v)|\Phi_u(T_1')|,\\
F(T_2;1,1)&=&r+F(T_1';f_2,1)=r+2^r+F(T_1'-u;1,1)+2^r|\Phi_u(T_1')|.
\end{array}
$$

Hence we have
$$F(T_2;1,1)-F(T_1;1,1)=[2^r+r-F(T_2';1,1)]+[2^r-F(T_2';1,1;v)]|\Phi_u(T_1')|.$$
Note that $T_2'$ is a tree with $r+1$ vertices. Hence, by Theorem
1.1 or Corollary 2.8,
$$2^r+r-F(T_2';1,1)\geq 0$$
with equality holds if and only if $T_2'=K_{1,r}$. Since $T_2'$
has at least $r$ subtrees $v_i$'s ($v_i\neq v$) with a vertex,
each of which is not a subtree of $T_2'$ containing vertex $v$,
$$F(T_2';1,1)\geq F(T_2';1,1;v)+r.$$Note that $F(K_{1,r};1,1)=2^r+r.$
Hence
$$0\leq F(K_{1,r};1,1)-F(T_2';1,1)\leq 2^r-F(T_2';1,1;v).$$
Therefore, we have
$$2^r-F(T_2';1,1;v)\geq 0$$
with equality holds if and only if $T_2'=K_{1,r}$ and
$d_{T_2'}(v)=r$. Hence we have
$$F(T_2;1)\geq F(T_1;1)$$
with equality holds if and only if $T_2'=K_{1,r}$ and
$d_{T_2'}(v)=r$. The Lemma thus follows. $\blacksquare$

Let $T_1'$ and $T_2'$ be two trees, and let $u$ (resp. $v$) be a
vertex of $T_1'$ (resp. $T_2'$), where $|V(T_2')|=r+1\geq 2$. Let
$T_1$ be the tree defined as above (see Figure 5(a)). Construct a
tree $T_3$ from $T_1'$ by identifying vertex $u$ of $T_1'$ and one
of two pendant vertices of a path with $r+1$ vertices (see Figure
5(c)). We call the procedure constructing $T_3$ from $T_1$ the
second transformation of tree $T_1$, denoted by
$\phi_{2}(T_1)=T_3$. As that in the proof of Lemma 3.1 we can
prove the following:
\begin{lem}
Let $T_1$ and $T_3$ be the trees defined as above, where $r\geq
1$. Then
$$\chi(T_1)=F(T_1;1,1)\geq \chi(T_3)=F(T_3;1,1)$$
with equality holds if and only if $T_2'=P_{r+1}$ and
$d_{T_2'}(v)=1$.
\end{lem}
\begin{rem}
Let $T$ be a tree with $n$ vertices and $T\neq K_{1,n-1}$ and
$T\neq P_{n}$. Suppose that $(v',u)$ is a pendant edge of $T$ and
$d_T(v')=1$. Let $T_1'$ be the subtree of $T$ containing two
vertices $v'$ and $u$, and let $T_2'=T-v'$. Obviously, with
application of the first (resp. second) transformation of tree
$T$, $T$ can be transformed into the star $K_{1,n-1}$ (resp. the
path $P_n$). Hence, by Lemma 3.1 (resp. Lemma 3.2),
$F(T;1,1)<F(K_{1,n-1};1,1)$ (resp. $F(T;1,1)>F(P_n;1,1)$).
\end{rem}

\par Suppose $V(P_{d+1})=\{v_1,v_2,\ldots,v_{d+1}\}$ and
$E(P_{d+1})=\{(v_j,v_{j+1})|j=1,2,\ldots,d\}$ are the vertex set
and edge set of a path $P_{d+1}$ with $d+1$ vertices,
respectively. Assume that $k_i,k_{i+1},\ldots,k_d$ are $d-i+1$
non-negative integers and $k_i\neq 0$. Construct two trees,
denoted by $T=T_d(k_i,k_{i+1},\ldots,k_d)$ and
$T^*=T_d(k_i+k_{i+1},k_{i+2},\ldots,k_d)$, with
$d+1+\sum\limits_{l=i}^dk_l$ vertices as follows. $T$ is the tree
obtained from $P_{d+1}$ by attaching $k_l$ pendant edges to
vertices $v_l$ for $l=i,i+1,\ldots,d$ (see Figure 6(a)) and $T^*$
is the tree obtained from $P_{d+1}$ by attaching $k_i+k_{i+1}$
pendant edges to vertex $v_{i+1}$ and $k_l$ pendant edges to
vertices $v_l$ for $l=i+2,i+3,\ldots,d$ (see Figure 6(b)). We call
the procedure constructing $T^*$ from $T$ the third transformation
of tree $T$, denoted by $\phi_3(T)=T^*$.
\begin{figure}[htbp]
  \centering
 \scalebox{1}{\includegraphics{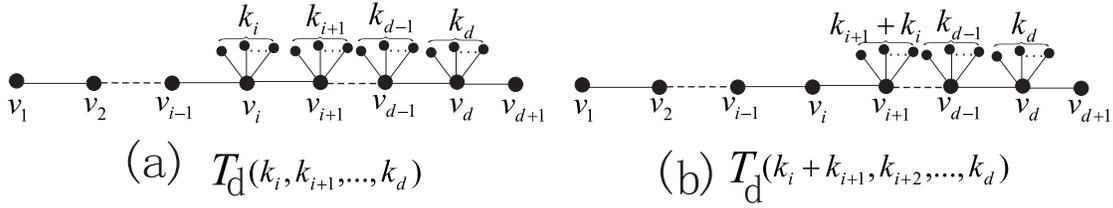}}
  \caption{\ \ (a)\ The tree
$T_d(k_i,k_{i+1},\ldots,k_d)$. \ (b)\ The tree
$T_d(k_i+k_{i+1},k_{i+2},\ldots,k_d)$.}
\end{figure}
\begin{lem}
Suppose $d$ and $k_l$ for $l=i,i+1,\ldots,d$ are non-negative
integers and $d>1, k_i\geq 1$. Let $T=T_d(k_i,k_{i+1},\ldots,k_d)$
and $T^*=T_d(k_i+k_{i+1},k_{i+2},\ldots,k_d)$ be the two trees
defined as above. If $i\leq \frac{d+1}{2}$, then we have
$$F(T;1,1)\leq F(T^*;1,1)$$
with equality holds if and only if $k_{i+1}=k_{i+2}=\ldots=k_d=0$,
$d$ is odd and $i=\frac{d+1}{2}$.
\end{lem}

{\bf Proof}\ \ We assume that $T_1$ is one of two components of
$T-(v_{i+1},v_{i+2})$, which contains vertex $v_{d+1}$. Obviously,
$T_1$ is a subtree of $T$ and it can be naturally regarded as a
subtree of $T^*$. By Algorithms 2.9 and 2.10, we have
$$
\begin{array}{lll}
F(T;1,1)&=&\frac{1}{2}(i-1)i+k_i+k_{i+1}+F(T_1;1,1)+i2^{k_i}+2^{k_{i+1}}[F(T_1;1,1;v_{i+2})+1]\\
&&+i2^{k_i+k_{i+1}}[F(T_1;1,1;v_{i+2})+1];\\
F(T^*;1,1)&=&\frac{1}{2}i(i+1)+k_{i}+k_{i+1}+F(T_1;1,1)+2^{k_i+k_{i+1}}(i+1)[F(T_1;1,1;v_{i+2})+1].
\end{array}
$$
Hence it is easy to obtain the following
$$F(T^*;1,1)-F(T;1,1)=[2^{k_i}-1][2^{k_{i+1}}F(T_1;1,1;v_{i+2})+2^{k_{i+1}}-i].$$
Note that $k_i>0$. So we have $2^{k_i}-1>0$. Since
$F(T_1;1,1;v_{i+2})$ has at least $d+1-(i+1)=d-i$ vertices,
$F(T_1;1,1;v_{i+2})\geq d-i$, which implies that
$$2^{k_{i+1}}F(T_1;1,1;v_{i+2})+2^{k_{i+1}}-i\geq 2^{k_{i+1}}(d-i+1)-i\geq d-2i+1$$
with equality holds if and only if $k_{i+1}=0$ and
$F(T_1;1,1;v_{i+2})=d-i$. Since $i\leq \frac{d+1}{2}$, we have
$$2^{k_{i+1}}F(T_1;1,1;v_{i+2})+2^{k_{i+1}}-i\geq
2^{k_{i+1}}(d-i+1)-i\geq d-2i+1\geq 0$$ with equality if and only
$k_{i+1}=0, i=\frac{d+1}{2}$ and $F(T_1;1,1;v_{i+2})=d-i$. It is
not difficult to see $F(T_1;1,1;v_{i+2})=d-i$ if and only if
$k_{i+1}=k_{i+2}=\ldots=k_d=0$. Hence we have prove that
$F(T;1,1)\leq F(T^*;1,1)$ with equality holds if and only if
$k_{i+1}=k_{i+2}=\ldots=k_d=0$, $d$ is odd and $i=\frac{d+1}{2}$.
Hence the lemma follows. $\blacksquare$

Let $T_0$ be a tree with at least two vertices and $u$ a vertex of
$T$. For arbitrary two positive integers $s,t$, construct a tree,
denoted by $T_0(s,t)$, from $T_0$ by attaching two paths with
$s+1$ and $t+1$ vertices to vertex $u$. Figure 7(a) and (b)
illustrate two trees $T_0(s,t)$ and $T_0(s+t-1,1)$. We call the
procedure constructing $T_0(s+t-1,1)$ from $T_0(s,t)$ the fourth
transformation of $T_0(s,t)$, denoted by
$\phi_4(T_0(s,t))=T_0(s+t-1,1)$.
\begin{figure}[htbp]
  \centering
 \scalebox{1}{\includegraphics{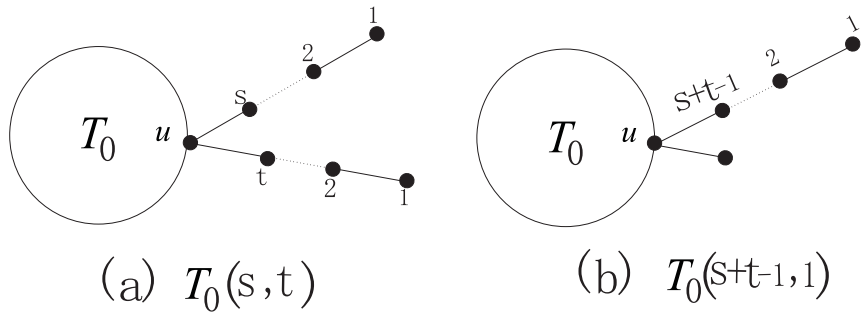}}
  \caption{\ \ (a)\ The tree $T_0(s,t)$. \ (b)\ The
tree $T_0(s+t-1,1)$.}
\end{figure}
\begin{lem}
Let $T_0$ be a tree with at least two vertices and $u$ a vertex of
$T_0$. For arbitrary two positive integers $s\geq 2,t\geq 2$, let
$T_0(s,t)$ be the tree defined as above. Then
$$F(T_0(s,t);1,1)>F(T_0(s+t-1,1);1,1).$$
\end{lem}

{\bf Proof}\ \ Let $f_i: V(T_0)\longrightarrow \mathcal R$
($i=1,2$) be two functions defined as follows:
$$
f_1(v)=\left\{
\begin{array}{ll}
(s+1)(t+1) & \ \mbox{if}\ \ v=u,\\
1 & \ \mbox{otherwise}.
\end{array}
\right. , f_2(v)=\left\{
\begin{array}{ll}
2(s+t) & \ \mbox{if}\ \ v=u,\\
1 & \ \mbox{otherwise}.
\end{array}
\right..
$$
Suppose that $\Phi_u(T_0)$ is the set of subtrees of $T_0$ with as
least two vertices, each of which contains vertex $u$. By
Algorithms 2.9 and 2.10, we have
$$
F(T_0(s,t);1,1)=\frac{1}{2}s(s+1)+\frac{1}{2}t(t+1)+F(T_0-u;1,1)+(s+1)(t+1)+(s+1)(t+1)|\Phi_u(T_0)|;
$$
$$
F(T_0(s+t-1,1);1,1)=1+\frac{1}{2}(s+t-1)(s+t)+F(T_0-u;1,1)+2(s+t)+2(s+t)|\Phi_u(T_0)|.
$$
From the equalities above, we have
$$F(T_0(s,t);1,1)-F(T_0(s+t-1,1);1,1)=(st-s-t+1)|\Phi_u(T_0)|.$$
Since $s\geq 2$ and $t\geq 2$, we have $st>s+t-1$. Hence
$$(st-s-t+1)|\Phi_u(T_0)|>0$$ which implies
$$F(T_0(s,t);1,1)>F(T_0(s+t-1,1);1,1).$$
Hence we have finished the proof of the lemma. $\blacksquare$
\subsection{Trees with extremal number of subtrees}
\hspace*{\parindent}
First, we need to defined two trees as follows. Suppose $n,d$ and
$\Delta$ are three positive integers, $n\geq d+1$ and $\Delta\geq
2$. Let $T_{n,\Delta}$ be the tree obtained from $P_{n-\Delta+1}$
by attaching $\Delta-1$ pendant edges to one of pendant vertices
of $P_{\Delta-1}$ (see Figure 8(a)).  Suppose
$V(P_{d+1})=\{1,2,\ldots,d+1\}$ and
$E(P_{d+1})=\{(i,i+1)|i=1,2,\ldots,d\}$ are the vertex set and
edge set of a path $P_{d+1}$ with $d+1$ vertices, respectively.
Let $T(n,d)$ be the tree obtained from $P_{d+1}$ by attaching
$n-d-1$ pendant edges to vertex $[\frac{d+1}{2}]+1$, where $[x]$
denotes the largest integer no more than $x$ (see Figure 8(b)).
\begin{figure}[htbp]
  \centering
 \scalebox{1}{\includegraphics{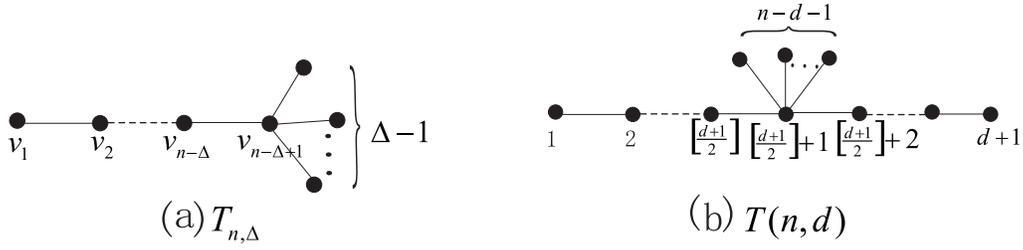}}
  \caption{\ \ (a)\ The tree $T_{n,\Delta}$. \ (b)\ The tree
  $T(n,d)$.}
\end{figure}
\begin{thm}
Let $\Delta$ be a positive integer more than two, and let $T$ be a
tree with $n$ vertices, which has the maximum degree at least
$\Delta$. Then
$$F(T;1,1)\geq F(T_{n,\Delta};1,1)$$
with equality holds if and only if $T=T_{n,\Delta}$, where
$T_{n,\Delta}$ is the tree defined as above.
\end{thm}
\begin{thm}
Let $d$ be a positive integer more than one, and let $T$ be a tree
with $n$ vertices, which has diameter at least $d$. If $T\neq
T(n,d)$, then
$$F(T;1,1)<F(T(n,d);1,1),$$
where $T(n,d)$ is the tree defined as above.
\end{thm}

Before we prove the theorems above, we consider some of their
corollaries, which characterize the trees with the second, third,
fourth, and fifth largest numbers of subtrees and the tree with
the second minimum number of subtrees.

Since the maximum degree of a tree $T$ with $n$ vertices, which is
different from $P_n$, is more than two, the following corollary is
immediate from Theorems 3.6 and 1.1.
\begin{cor}
Let $T$ be a tree with $n$ ($n\geq 3$) vertices and $T\neq P_n,
T\neq T_{n,3}$. Then
$$F(T;1,1)>F(T_{n,3};1,1)>F(P_n;1,1).$$
\end{cor}

In order to present Corollary 3.9, we need to define a new tree
$B_{n,d}$ (where $n\geq 2d+2\geq 4$) as follows. Let $B_{n,d}$ be
the tree with $n$ vertices obtained from $K_{1,n-d-1}$ by
attaching $d$ pendant edges to one of pendant vertices of
$K_{1,n-d-1}$ (Figures 9(b) and (c) show $B_{n,2}$ and $B_{n,3}$,
respectively). Obviously, $B_{n,1}=T(n,3)$ (see Figure 9(a)).
\begin{figure}[htbp]
  \centering
 \scalebox{1}{\includegraphics{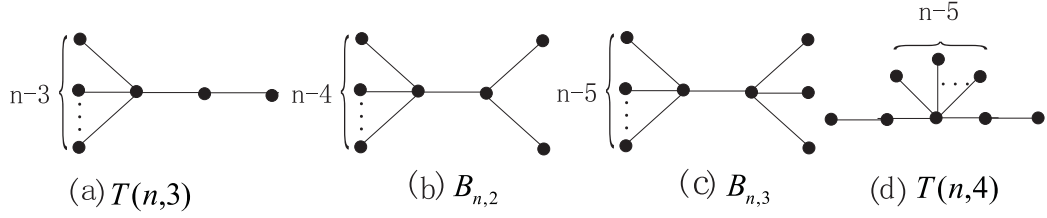}}
  \caption{\ \ (a)\ The tree $B_{n,1}=T(n,3)$. \ (b)\ The tree
$B_{n,2}$. \ (c)\ The tree $B_{n,3}$.\ (d)\ The tree $T(n,4)$.}
\end{figure}
\begin{cor}
Let $T$ be a tree with $n\geq 8$ vertices and $T\neq
K_{1,n-1},T(n,3),B_{n,2},B_{n,3}$,$T(n,4)$ (see Figure
9(a)$-$(d)). Then
$$F(K_{1,n-1};1,1)>F(T(n,3);1,1)>F(B_{n,2};1,1)$$$$>F(B_{n,3};1,1)>F(T(n,4);1,1)>F(T;1,1).$$
\end{cor}

{\bf Proof}\ \  By Theorem 1.1 and Theorem 3.7, we have
$$F(K_{1,n-1};1,1)>F(T(n,3);1,1)>F(B_{n,2};1,1). \eqno{(13)}$$
If the diameter of $T$ is at least $4$, then by Theorem 3.7 we
have
$$F(T(n,4);1,1)>F(T;1,1). \eqno{(14)}$$
The following equalities can be proved from Algorithm 2.9:
$$F(B_{n,d};1,1)=n-2+2^d+2^{n-d-2}+2^{n-2}, \eqno{(15)}$$
$$F(T(n,4);1,1)=n+1+2^{n-2}+2^{n-5}. \eqno{(16)}$$ Obviously, if $n\geq 8$, then by
$(15)$ and $(16)$ we have
$$F(B_{n,2};1,1)>F(B_{n,3};1,1),\ F(B_{n,3};1,1)>F(T(n,4);1,1).$$

Note that if the diameter of a tree $T'$ with $n=8$ or $n=9$
vertices equals three, then $T$ must be one of $K_{1,n-1}, T(n,3),
B_{n,2}$, and $B_{n,3}$. Hence the corollary holds when $n=8$ or
$n=9$.

Note that if the diameter of a tree $T'$ with $n\geq 10$ vertices
equals three, then $T'$ must has the form of $B_{n,d}$, where
$n\geq 2d+2$ (by the definition of $B_{n,d}$). By $(15)$ and
$(16)$,
$$F(B_{n,i};1,1)-F(T(n,4);1,1)=2^i+2^{n-i-2}-3-2^{n-5}.$$
By the definition of $B_{n,i}$, $n\geq 2i+2$. It is not difficult
to show that if $n\geq 2i+2\geq 10$ (hence $i\geq 4$), then
$$F(B_{n,i};1,1)<F(T(n,4);1,1).$$
Therefore, we have shown that if $n\geq 10$ and $i\geq 4$, then
$$F(B_{n,1};1,1)>F(B_{n,2};1,1)>F(B_{n,3};1,1)>F(T(n,4);1,1)>F(B_{n,i};1,1). \eqno{(17)}$$

Hence the corollary follows. $\blacksquare$

{\bf Proof of Theorem 3.6}\ \ Let $T$ be a tree with $n$ vertices
and $T\neq T_{n,\Delta}$. Note that $T$ is a tree with the maximum
degree at least $\Delta$. Hence there exists a vertex $u$ of $T$
such that $d_T(u)\geq \Delta$. Without loss of generality, we
assume that $\{v_1,v_2,\ldots,v_{\Delta-1}\}$ is a subset of the
neighbor set of $u$ in $T$. Obviously, if we delete $\Delta-1$
edges $(u,v_1), (u,v_2),\ldots,(u,v_{\Delta-1})$ from $T$, then
$\Delta$ components $T_i$'s (for $i=1,2,\ldots,\Delta$) of $T$ can
be obtained, where $T_i$ is the component containing vertex $v_i$
for $i\leq \Delta-1$ and $T_{\Delta}$ is the one containing vertex
$u$. Furthermore, $T_{\Delta}$ contains at least two vertices.
Hence $T$ has the form illustrated in Figure 10(a).
\begin{figure}[htbp]
  \centering
 \scalebox{1}{\includegraphics{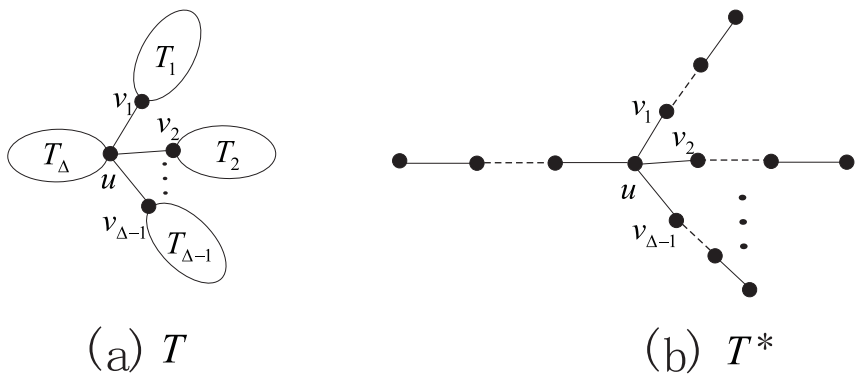}}
  \caption{\ \ (a)\ The tree $T$ in the proof of
Theorem 3.6. \ (b)\ The tree $T^*$ in the proof of Theorem 3.6.}
\end{figure}


With repeated applications of the second transformations of trees,
$T$ can be transformed to the form of $T^*$ showed in Figure
10(b). Hence by Lemma 3.2 we have $F(T;1,1)>F(T^*;1,1)$. If
$T^*=T_{n,\Delta}$, then the theorem holds. If $T^*\neq
T_{n,\Delta}$, then by repeated applications of the forth
transformations of trees $T^*$ can be transformed to
$T_{n,\Delta}$, and we have $F(T^*;1,1)>F(T_{n,\Delta};1,1)$.
Hence $F(T;1,1)>F(T_{n,\Delta},1,1)$. The theorem thus has been
proved. $\blacksquare$

{\bf Proof of Theorem 3.7}\ \ Let $T$ be a tree with $n$ vertices
with the diameter at least $d$ and $T\neq T(n,d)$. Then there
exists a path of length $d-1$ in $T$, denoted by
$P_{d}=P(v_1-v_2-\ldots-v_d)$, where $d_T(v_1)=1$. Then $T$ must
has the form illustrated in Figure 11(a), where $T_i$ is a subtree
of $T$ containing vertex $v_i$ for $i=2,3,\ldots,d$. Particularly,
since the diameter of $T$ is at least $d$, $T_d$ contains at least
two vertices. With repeated applications of the first
transformations of trees, $T$ can be transformed to the tree with
form of $T^*$ shown in Figure 11(b) and hence we have the
following:
$$F(T;1,1)\leq F(T^*;1,1)$$
with equality holds if and only if $T=T^*$.
\begin{figure}[htbp]
  \centering
 \scalebox{1}{\includegraphics{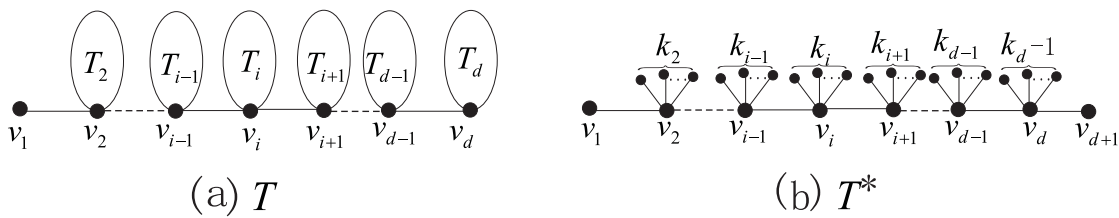}}
  \caption{\ \ (a)\ The tree $T$ in the proof of
Theorem 3.7. \ (b)\ The tree $T^*$ in the proof of Theorem 3.7.}
\end{figure}

If $T^*\neq T(n,d)$, then by repeated applications of the third
transformations of trees $T^*$ can be transformed to $T(n,d)$ and
hence $F(T^*;1,1)<F(T(n,d);1,1)$. So $F(T;1,1)<F(T(n,d);1,1)$. If
$T^*=T(n,d)$, then $T\neq T^*$. But in this case we have shown
that $F(T;1,1)<F(T^*;1,1)=F(T(n,d);1,1)$. The theorem thus
follows. $\blacksquare$
\section{Concluding remarks}
\hspace*{\parindent}
In this paper, we have investigated the problem on enumeration of
subtrees of trees. We obtained a linear-time algorithm to count
the sum of weights of subtrees of a tree and we also characterized
some trees with extremal number of subtrees. Note that if $G$ is a
connected graph then some coefficients of its Tutte polynomial
$T_G(x,y)$ can count the numbers of some kinds of subgraphs of $G$
\cite{Boll98}. For example, $T_G(1,1)$ is the number of spanning
trees of $G$, $T_G(2,1)$ is the number of forests in $G$,
$T_G(1,2)$ is the number of connected spanning subgraphs in $G$,
and $T_G(2,2)$ equals the number of spanning subgraphs in $G$. A
natural extension of our work would be to give some methods to
enumerate connected subgraphs of a connected graph. On the other
hand, an acyclic molecular can be expressed by a tree in quantum
chemistry (see \cite{GC89}). The study of the topological indices
(see, for example, \cite{KMPT92,RZ01}) has been undergoing rapid
expansion in the last few years. Obviously, the number of subtrees
of a tree can be regarded as a topological index. Hence another
interesting direction is to explore the role of this index in
quantum chemistry.
\newpage
\vskip0.5cm \noindent {\bf Acknowledgements}

We are grateful to the anonymous referee for the comments that
greatly improved the presentation of the paper.


\end{document}